\documentclass[a4paper]{paper}

\usepackage[utf8]{inputenc}
\usepackage[english]{babel}
\usepackage{amsmath}
\usepackage{amsthm}
\usepackage{amssymb}
\usepackage{amsfonts}
\usepackage{graphicx} 
\usepackage{csquotes}
\usepackage{biblatex}
\usepackage[nottoc,numbib]{tocbibind}
\usepackage{hyperref}
\hypersetup{
    colorlinks=true,
    linkcolor=blue,
    filecolor=red,      
    urlcolor=blue,
    citecolor=red,
    } 
\usepackage[nameinlink]{cleveref} % Must be loaded after hyperref
\usepackage{acro}
\usepackage{floatpag}
\usepackage{geometry}
\usepackage{todonotes}
\usepackage{attachfile}
\usepackage{blindtext}
\usepackage{multicol}

\newcommand{\Rwidth}{0.21}

\newcommand{\real}{\mathbb{R}}

\newcommand{\Natural}{\mathbb{N}}
\newcommand{\DataManifold}{M}

\newcommand{\matrixForm}[1]{\mathsf{#1}}

\newcommand{\Fourier}{\operatorname{\mathcal{F}}}
\newcommand{\Radon}{\operatorname{\mathcal{P}}}
\newcommand{\ReLU}{\operatorname{ReLU}}
\newcommand{\FNO}{\operatorname{FNO}}
\def\imagetop#1{\vtop{\null\hbox{#1}}}

\DeclareAcronym{CT}{
  short = CT,
  long = computed tomography}
\DeclareAcronym{FBP}{
  short = FBP,
  long = filtered back projection}
\DeclareAcronym{FDK}{
  short = FDK,
  long = {Feldkamp–Davis–Kress}}
\DeclareAcronym{FNO}{
  short = FNO,
  long = Fourier neural operator}
\DeclareAcronym{PSNR}{
  short = PSNR,
  long = peak signal-to-noise ratio}
\DeclareAcronym{SSIM}{
  short = SSIM,
  long = structural similarity index}  
\DeclareAcronym{GPU}{
  short = GPU,
  long = graphics processing unit}  
\DeclareAcronym{FNOBP}{
  short = FNO-BP,
  long = FNO back projection} 
\DeclareAcronym{GELU}{
  short = GELU,
  long = Gaussian error linear unit}

\title{Fast deep learning based  reconstruction for limited angle tomography}

\author{Knut Salomonsson \and Eric Oldgren \and Emanuel Ström \and Ozan Öktem}

\begin{document}
\maketitle

\begin{abstract}
A major challenge in computed tomography is reconstructing objects from incomplete data. An increasingly popular solution for these problems is to incorporate deep learning models into reconstruction algorithms. This study introduces a novel approach by integrating a \ac{FNO} into the Filtered Backprojection (FBP) reconstruction method, yielding the \ac{FNOBP} network. We employ moment conditions for sinogram extrapolation to assist the model in mitigating artefacts from limited data. Notably, our deep learning architecture maintains a runtime comparable to classical \ac{FBP} reconstructions, ensuring swift performance during both inference and training. We assess our reconstruction method in the context of the Helsinki Tomography Challenge 2022 and also compare it against regular \ac{FBP} methods.

\end{abstract}
\acresetall

\section{Background, motivation, and specific contributions}
\subsection{Principle of tomographic imaging}
Tomography is a technique to non-invasively image the interior structure of an object by repeatedly probing the object with particles/waves from directions that surround it. 
These particles/waves are detected after they have interacted with the (interior structure) object. Observed data, henceforth called a sinogram, therefore represents noisy indirect observations of the interior structure one seeks to determine. Tomographic imaging therefore includes a reconstruction method, which is a computational method for recovering the interior structure from a sinogram.

A well-known example of the above is x-ray \ac{CT} where one probes an object with x-ray photons. The nature of the object depends on the use case, e.g., in medical imaging, it is a human being and in non-destructive testing, it is commonly a manufactured object.
The interior structure is represented by the (spatially varying) linear attenuation coefficient, which is commonly visualized as an \emph{image}.
Data (\emph{sinogram}) consists of x-ray intensity measurements along a specific set of lines passing through the object. The \emph{acquisition geometry} refers to a mathematical specification of the arrangement of these lines (manifold of lines).

\subsection{Evolution of \ac{CT} reconstruction methods}
The reconstruction method used for recovering the image from a sinogram is a key part of \ac{CT} imaging. 
Many reconstruction methods are derived in a continuum setting where the image one seeks to recover is represented by a real-valued function on 2- or 3-dimensional space and the sinogram represents samples of a real-valued function on lines. 
Furthermore, in \ac{CT} one commonly adopts the Beer--Lambert law for modeling how X-ray photons interact with tissue/matter. Combining this simplified physics model with appropriate pre-processing of measured data allows one to re-phrase the reconstruction task as inverting the ray transform restricted to the manifold of lines given by the acquisition geometry.

We next describe the evolution of \ac{CT} reconstruction methods that are based on the above reformulation of the reconstruction task.
From a mathematical viewpoint, the first steps towards reconstruction were taken already in 1917 when Johan Radon published an analytic inversion formula for the Radon transform \cite{Radon:1917aa}, which gives an inversion formula for the ray transform in the planar setting. 
A similar explicit inversion formula for the ray transform in general dimensions was derived in 1938 by Fritz John \cite{John:1938aa}. 
Since then, much work in the mathematics community has been devoted to characterizing injectivity, range, and stability of such integral transforms as well as deriving explicit inversion formulas for them, see \cite{Markoe:2006aa} for a survey of this mathematical development and \cite{Webb1990} for en extensive account of the historical development. 
% The work by Johann Radon and other mathematicians was to large extent unknown to practitioners involved in the development of tomography in the 1950-1960s. 

An exact inversion formula for the ray transform (when it is restricted to lines from a specific acquisition geometry) allows one to exactly recover an image from a sinogram that represents noise-free continuum data. 
Inverting the ray transform is however an \emph{ill-posed} problem, which means that exact inversion will inevitably involve steps that amplify small measurement/sampling errors in the sinogram.
Exact inversion formulas have therefore had limited practical applicability, but they serve as a basis for \emph{analytic (reconstruction) methods}. These methods use principles from Fourier analysis and sampling theory to stabilize the intrinsically unstable inversion, e.g., by recovering the band-limited part of the image whereas high-frequency components, like noise, are filtered out.
Such methods were pioneered by Bracewell and Riddle, who in 1967 published the \ac{FBP} method for parallel lines in the planar setting \cite{Bracewell:1967aa}.
Since then, many variants of \ac{FBP} have been developed, mainly for adapting it to other acquisition geometries \cite{Shepp:1974aa,Natterer:1980aa, Feldkamp:1984aa,Faber:1995aa,Katsevich:2003aa,Noo:2003aa,Katsevich:2006aa}.
Another closely related class of analytic methods are Fourier methods which represent a direct implementation of the projection slice theorem (which relates the ray transform to the Fourier transform). 
 
In addition to analytic methods, there are three additional classes of reconstruction algorithms for inverting the ray transform. One is \emph{model-based iterative methods}, like algebraic reconstruction technique, simultaneous iterative reconstruction technique, and the least-squares conjugate gradient method. Here, the ray transform is approximately inverted by an iterative scheme, which in the limit is designed to converge to an image that maximizes consistency against the sinogram. The typical setup is to use an iterative scheme that minimizes a data fidelity term, as in computing a least-squares solution. Iterates are however stopped long before convergence (early stopping) in order to avoid instabilities that would otherwise arise due to ill-posedness, see \cite[sec.~5.3]{NattererWubbeling2001} for more details.
Another class is the \emph{model-based variational methods}, like quadratic Tikhonov or total variation regularisation. Here, one minimizes an objective that combines a data fidelity term with a regulariser. The data fidelity ensures consistency against data (as in model-based iterative methods) whereas the regulariser penalizes instabilities, see 
\cite{OtmarBook,BurgerReview} for an extensive survey of these methods. 
The final class of reconstruction methods are those that are based on \emph{deep learning}. The starting point is to view a reconstruction method as a statistical estimator in Bayesian inference.
The estimator, which is a reconstruction method (not a reconstruction) is then learned from training data. 
The type of training data one has access to and how one chooses to set up the learning problem determines which estimator one seeks to learn. The deep neural network architecture is here a parametrization of possible estimators. 
It is common to use domain-adapted architectures that encode some of the underlying physics since these come with vastly better generalization properties than generic deep neural network architectures, see \cite[sec.~5]{Arridge:2019aa} for more on these matters.

\subsection{Challenges for clinical \ac{CT} reconstruction}
The overall aim in \ac{CT} reconstruction is to have a reconstruction method that produces `high quality' 3D images in a timely manner and with a memory footprint that does not require using supercomputing resources.

For full-dose 3D clinical \ac{CT}, analytic methods based on \ac{FBP} offer a compromise in balancing image quality against reconstruction speed and memory footprint that is still difficult to outperform \cite{Pan:2009aa}. 
However, a drawback with analytic methods is that they need to be specifically tailored to the data acquisition geometry, e.g., it is difficult to make use of \emph{all} the measured data without resorting to approximations.
The image quality obtained from analytic methods becomes an issue when \ac{CT} data is either under-sampled (sparse view \ac{CT}) or unevenly sampled, as in limited angle \ac{CT}.
The stabilization strategy used in analytic methods, which is sufficient for full-dose \ac{CT}, is simply not strong enough to provide useful reconstructions in these settings. 
This is especially notable in limited angle \ac{CT} where the x-ray source in the \ac{CT} scanner cannot move all the way around the object. This results in a sinogram that has values of the ray transform on lines only from a limited range of directions, a setting that is known to significantly increase the ill-posedness of the reconstruction problem.
Another final drawback of analytic methods is that they are inherently deterministic, so they are not suitable to account for statistical properties in the noise. This is an issue when data is highly noisy (low-dose \ac{CT}).

Attempts to address these drawbacks of analytic methods have catalyzed the development of model-based iterative/variational methods, some of which have made it into clinic practice as surveyed in \cite{Willemink:2019aa,LaRivierea:2021aa}.
Compared to analytic methods, iterative/variational methods offer the possibility to incorporate more sophisticated (and stronger) stabilization strategies. 
Next, unlike analytic methods, they are also not tailored for a specific acquisition geometry.
Furthermore, they can incorporate a more accurate physics model for the interaction of X-ray photons with tissue/matter than the one provided by the Beer--Lambert law that is used in analytic methods. Finally, iterative/variational methods can also account for statistical properties of the noise in data, which is essential in low-dose \ac{CT}.
On the other hand, these model-based methods have a long run-time that becomes prohibitive in time-critical use cases, like \ac{CT} imaging in trauma. 

Deep learning-based methods form a final category of reconstruction.
State-of-the-art approaches that rely on domain adapted physics informed deep neural network architectures \cite{Adler2018} outperform model-based methods regarding the image quality while having a reconstruction run-time that is closer to the analytic methods \cite{Arridge:2019aa,Willemink:2019aa,LaRivierea:2021aa}.
With clever engineering, one can limit their memory footprint and thereby train these domain-adapted architectures on high-end gaming \acp{GPU} hardware also in the 3D setting \cite{LDP3D}.
This does however not address the time that the training takes. The main bottleneck here is that the training involves computing a large amount of forward projections (evaluations of the ray transform on all lines in the sinogram) and corresponding back projections. 

This brings us to the \emph{challenge} considered in this paper, which is to develop a \emph{deep neural network-based reconstruction method that provides a significantly better reconstruction quality than analytic methods (in limited data setting) while having a run-time and memory footprint that is comparable to analytic methods}. 
This requirement constrains the deep neural network for reconstruction, e.g., the highly successful architectures for reconstruction derived by unrolling some iterative scheme \cite{Adler2018,Adler:2023aa,Arridge:2019aa} involves several forward- and back-projections, so these are for this reason disqualified. 
In fact, having a run-time and memory footprint that is comparable to analytic methods means using an architecture that avoids forward projection(s) and incorporates a single (or a few) back projection(s). 

\subsection{Specific contributions}
One contribution is sinogram extrapolation in fan-beam limited angle tomography which is based on the range conditions for the ray transform. 
Another contribution is to introduce the novel \ac{FNOBP} network architecture that is specifically adapted for tomographic reconstruction. It avoids forward projections and incorporates a single back projection, so its run-time and memory footprint during inference are comparable to \ac{FBP}. Meanwhile, it outperforms \ac{FBP} regarding reconstruction quality when tested on the HTC-2022 limited angle tomography data \cite{HTC2022_data}.

\section{Computationally efficient reconstruction methods} \label{method}

We here provide a brief overview of computationally efficient reconstruction methods. This refers to methods for recovering an image from a sinogram that only involves computing a single (or a few) back-projection operation. In addition, one avoids computing forward projections altogether. 
This stands in contrast to model-based iterative and variational methods and most of the state-of-the-art deep learning-based methods, which rely on domain-adapted architectures that incorporate forward projections \cite{Adler:2023aa,Arridge:2019aa}.

\subsection{Beer-Lambert law}
When radiating an object with x-ray photons, the photons will interact with the spatially varying linear attenuation coefficient of the object. This is the `interior structure' of the object that one can sense with x-rays and it is mathematically represented by a function $f \colon \real^d \to \real$ with $d=2$ if all x-rays are in a fixed planar 2D cross-section and $d=3$ otherwise. 

The radiative transport equation yields an accurate model for the interaction between the object and x-ray photons that travel along a line $\ell$ through the object \cite{Bal:2009aa}.
Simplifying this model (assuming monochromatic x-rays and disregarding scattering events) results in the Beer-Lambert law, which states that  
\begin{equation}\label{eq:TomoData} 
  \frac{I_{\text{out}}(\ell)}{I_{\text{in}}(\ell)} =
    e^{\int_{\ell} f(x)\,d\mu(x)}
  \implies 
  g(\ell) := -\log\biggl( \frac{I_{\text{out}}(\ell)}{I_{\text{in}}(\ell)} \biggr) = \int_{\ell} f(x)\,d\mu(x)
  \quad\text{for $\ell \in \DataManifold$.}
\end{equation}
In the above, $\DataManifold$ is a collection of lines through the object, $\mu$ is the 1D Lebesgue measure on a line $\ell \in \DataManifold$, $I_{\text{in}}(\ell)$ is the intensity of x-ray photons traveling along $\ell$ before they interact with the object, and $I_{\text{out}}$ is the corresponding intensity x-ray photons after they have interacted with the object.
In general, during a tomographic data acquisition one measures $I_{\text{out}}(\ell)$ (transmission) along finitely many  lines $\ell \in \DataManifold$ and the corresponding incident intensity $I_{\text{in}}(\ell)$. The data acquisition geometry specifies the set $\DataManifold$ along with the actual sampling of lines in $\DataManifold$ that corresponds to measured data. The term $g(\ell)$, which is computed from data, is called a projection (or absorption) along $\ell$. A sinogram is the function $g \colon \DataManifold \to \real$, i.e., the collection of all projections as $\ell \in \DataManifold$ varies.

\subsection{Notions from integral geometry}
Describing the deep neural network architectures considered in this paper requires us to introduce some notions from integral geometry. 

We start with defining the \emph{ray-transform} restricted to a given set of lines $\DataManifold$ as 
\[
\Radon(f)(\ell) := 
  \int_{\ell} f(x)\,d\mu(x)
  \quad\text{for a line $\ell \in \DataManifold$.}
\]
The ray transform is mathematically a linear operator that maps a real-valued function $f \colon \real^d \to \real$ to a real-valued function on lines in $\DataManifold$. It is well-defined whenever $f$ is measurable and has compact support (it is in fact enough to require that $f$ decays faster than any polynomial, see \cite[sec.~2.1]{NattererWubbeling2001}).
Next, the Beer--Lambert law allows us to view the reconstruction task in \ac{CT} as inverting the ray transform restricted to $\DataManifold$. 
To see this, note that \ac{CT} reconstruction aims to recover the image $f \colon \real^d \to \real$ from the sinogram $g \colon \DataManifold \to \real$ in \eqref{eq:TomoData}. 
But $g(\ell) = \Radon(f)(\ell)$, so reconstruction is mathematically the same as inverting $\Radon$ restricted to $\DataManifold$.

The corresponding \emph{back-projection} is a key component in most reconstruction methods. 
When applied to a function $g \colon \DataManifold \to \real$, it is defined as   
\[
\Radon^*(g)(x) := 
  \int_{\DataManifold_x} g(\ell)\,d\sigma(\ell)
\quad\text{for $x \in \real^d$.}
\]
Here, $\DataManifold_x \subset \DataManifold$ denotes the subset of lines $\ell \in \DataManifold$ that pass through $x \in \real^d$ and $\sigma$ is the corresponding (Haar) measure on $\DataManifold_x$ (see \cite[sec.~3.3.4]{Markoe:2006aa}).
We see that $\Radon^*$ is a linear operator that maps functions $g \colon \DataManifold \to \real$ on lines onto real-valued functions on $\real^d$.
One can now show that $\Radon^*$ is the formal adjoint of the ray transform $\Radon$ (see \cite[theorem~3.29]{Markoe:2006aa}).

We finally mention the \emph{Radon transform}, which integrates $f \colon \real^d \to \real$ over hyper-planes (instead of lines) in $\real^d$. 
Since lines and hyper-planes coincide in $\real^2$, we get that the ray transform and the Radon transform coincide in the planar setting (when $d=2$). 

\subsection{\Acf{FBP}}\label{sec:FBP}
The \ac{FBP} method is based on the following identity \cite[theorem~2.13]{NattererWubbeling2001} that relates a  function $f$ to its sinogram $g := \Radon(f)$ on $\DataManifold$:
\[
  f \ast K = \Radon^*( k \circledast g ) 
  \quad
  \text{where $K := \Radon^*(k)$.}
\]
Here, $\ast$ in the left-hand-side is the $d$-dimensional convolution between functions on $\real^d$ whereas $\circledast$ in the right-hand-side is the $(d-1)$-dimensional convolution between functions on lines $\ell \in \DataManifold$. 
For the latter, the convolution evaluated at a line $\ell$ is given by a `common' $(d-1)$-dimensional convolution in the hyperplane through the origin that is orthogonal to $\ell$ \cite[eq.~(2.30)]{NattererWubbeling2001}.

Note that the sinogram $g=\Radon(f)$ is the noise-free continuum data, the latter meaning that the ray transform is assumed to be given on all lines in $\DataManifold$.
The function $k \colon \DataManifold \to \real$ is called the reconstruction kernel. 
It governs what feature of $f$ one seeks to recover. 
The \ac{FBP} method is said to be \emph{theoretically exact} whenever $k$ can be chosen so that $K = \delta$ since that corresponds to recovering $f$ exactly from continuum noise-free data. 
Inverting $\Radon$ is however an ill-posed problem, so computing the right-hand side will with such a choice of reconstruction kernel involve unstable steps. 
The idea in \ac{FBP} is to stabilize this by selecting a reconstruction kernel $k$ so that $K \approx \delta$.

Explicit reconstruction kernels exist for parallel and fan-beam geometries in planar settings. The same also holds for parallel beam geometries in 3D. The situation is more complicated for cone-beam geometries in 3D. Since the work by Katsevich \cite{Katsevich:2003aa}, there are explicit reconstruction kernels also in that setting with sources on any curve (circular, helical, \ldots). 
However, these theoretically exact \ac{FBP} methods do not make use of all the data, which in turn makes them overly sensitive to noise. 
For this reason, approximate \ac{FBP}-type of methods have been developed, like \ac{FDK} for cone-beam geometry with sources on a circle and variants of \ac{FDK} for other cone-beam acquisition geometries.

A specific challenge in the planar setting is when $\DataManifold$ is incomplete, which here means there is an open set on lines in $\real^2$ that are not in $\DataManifold$. 
In that setting, there is no good way of choosing a reconstruction kernel.
Applying standard \ac{FBP} reconstruction on a sinogram that is not given on an open subset means one implicitly assumes the sinogram is zero in that missing region. This is almost always inconsistent with measured data and thus results in very poor reconstructions. Sinogram extrapolation refers to methods for predicting the data in the missing region in order to pad with something other than zero. There are a few sinogram extrapolation methods, but these typically do not give adequate accuracy when they are combined with standard \ac{FBP}, so the resulting reconstructions are still quite poor. To improve upon this situation, one must combine the sinogram extrapolation with a modified reconstruction operator.

\subsection{Parallel and fan-beam coordinates}\label{sec:DataManifolds}
In a real-world situation, there is no access to the complete continuum data in $\DataManifold$. Instead, data is sampled according to a scheme that is determined by the acquisition geometry. As such one needs a parameterization of $M$ that is suitable for a given geometry. The most straightforward parametrization is the planar parallel beam geometry in which a line is parametrized by the angle of its direction and its distance from the origin. In such a geometry the Ray transform and its adjoint are explicitly given as
\begin{align*}
    \Radon (f)(\theta, s) &= \int_{-\infty}^\infty f\bigl( s\cos(\theta) - t\sin(\theta), s\sin(\theta) + t\cos(\theta) \bigr) dt \\
    \Radon^* (g) (x,y) &= \frac{1}{\pi} \int_0^\pi g\bigl( \theta, \cos(\theta)x + \sin(\theta)y \bigr) d\theta. 
\end{align*}
This parametrization determines the line $\ell$ as the line passing through the point $(s\cos(\theta), s\sin(\theta))$ perpendicular to the direction $(\cos(\theta), \sin(\theta))$.

While planar parallel beam geometry is the most simple to work with, it is not the most common. A more typical sampling scheme, and the one used when collecting dataset \cite{HTC2022_data,Meaney:2023aa}, is that of fan beam geometry with a flat detector and sources on a circular trajectory with radius $R>0$. In fan beam acquisition geometry, one coordinate determines the position of the line source on a circular orbit and another determines the direction of the line from that point. In this parametrization, the Ray transform can be written \cite{GuerreroBernardiMiqueles+2023+921+935}
\begin{multline*}
\Radon(f) (\beta, u) 
  = \int_{-\infty}^\infty f\bigl( R\cos(\beta)(1-t) + ut\sin(\beta), R\sin(\beta)(1-t) - ut\cos(\beta)\bigr) \sqrt{u^2+R^2} dt 
\\
\shoveleft{%
\Radon^* (g)(x,y) 
  = \frac{1}{2\pi} \int_0^{2\pi} g\Bigl( \beta, R\frac{x\sin\beta-y\cos\beta}{R-x\cos\beta-y\sin\beta} \Bigr) 
}
\\
\frac{R \sqrt{(x\sin\beta - y\cos\beta)^2 + (R - x\cos\beta - y\sin\beta)^2} }{(R - x\cos\beta - y\sin\beta)^2} d\beta.
\end{multline*}
 The two geometries can be related to one another via a change of coordinates. 
Let $(\theta, s)$ denote coordinates for lines in the parallel beam geometry. 
Likewise, let $(\beta, u)$ denote coordinates for lines in the fan beam geometry mentioned above.
We can now express the fan beam coordinates in terms of the parallel beam coordinates: 
\[
    \beta = \theta + \frac{\pi}{2} - \arctan \Bigl(\frac{s}{\sqrt{R^2 - s^2}}\Bigr) 
    \quad\text{and}\quad
    u = \frac{s R}{\sqrt{R^2-s^2}}.
\]
%The Jacobian for the change of variables from $(\theta, s)$ to $(\beta, u)$ is then 
%\[ \left|\frac{\partial (\theta, s)}{\partial (\beta, u)}\right| = \frac{R^3}{(u^2 + R^2)^{3/2}}.
%\]

\subsection{Sinogram extrapolation in the planar setting}\label{sec:sino_exp}

The aim is to extrapolate sinograms while accounting for the moment conditions that characterize the range of the ray transform. 
The method is an adaptation of \cite{Natterer2,HLCC_exp_chebyshev_lasso} to fan beam geometry (without re-binning).

The moment conditions for the planar ray transform are typically expressed in parallel beam coordinates \cite{Natterer2} as follows: 
A function $(\theta,s) \mapsto g(\theta, s)$ is in the range of $\Radon$ if and only if
\begin{equation}\label{eq:moment_conditions}
  %a_n(\theta) \stackrel{\text{def}}{=}
  \int_{-\infty}^{\infty} g(\theta, s) s^n ds \in \operatorname{span}
  \bigl\{e^{ik\theta} : 
      \text{$|k|\leq n, k+n$ even} 
  \bigr\}  
  \quad\text{for $n \in \Natural$.}
  %= \text{span}\{e^{-in\theta}, e^{-i(n-2)\theta}, ..., e^{in\theta} \} 
\end{equation}
To get a corresponding condition for other acquisition geometries, like fan beam, we note first that one can replace the polynomials $T_n(s) = s^n$ above with any set of polynomials $U_n(s)$ where $U_n$ is of degree $n$ and $U_n$ is even for $n$ even, and odd for $n$ odd. When $U_n$ are chosen to be an orthonormal family of polynomials with respect to some weight function $W(s)$ this imposes constraints on the series expansion of $g$ in the complete basis $U_n(s)W(s)e^{ik\theta}$. Hence, we get a new characterization of the range, see \cite{sino_exp_eul_lag, statistical_projection_completion} for details: $g$ lies in the range of $\Radon$ iff
\begin{equation}\label{eq:sino_synthesis}
    g(\theta, s) = \sum_{n}\sum_{k}c_{n,k}e^{ik\theta}U_n(s)W(s) 
\end{equation}
where $c_{n,k}$ is nonzero only if $|k|\leq n$ and $k + n$ is even. By neglected terms where $n>N$ we get an approximation of $g$ based on a finite set of coefficients $c = (c_{0,0}, c_{1,0}, \ldots, c_{N, N})^T$. In this form the moment conditions are easily transferred to other geometries, such as fan-beam.

We can exploit this to extrapolate a sub-sampled sinogram. If $g$ is known on a subset of its domain, we estimate coefficients $c_{n,k}$ to minimize the error in  \cref{eq:sino_synthesis} and then use the same coefficients to estimate $g$ in the rest of its domain. 
We can write this in matrix form in the discrete setting where $g$ is the array representing sinogram values:
\begin{equation*}
    g = \matrixForm{B} \cdot c.
\end{equation*}
Here, $\matrixForm{B}$ is a matrix whose columns consist of the basis functions $e^{ik\theta}U_n(s)W(s)$ sampled on a lattice in fan-beam coordinates (these values can be calculated using the coordinate mapping from  \cref{sec:DataManifolds}). 
Let $g_k$ be the sinogram samples from the known region and $g_u$ the sinogram in the unknown region.
\begin{equation*}
    g = \begin{bmatrix}
        g_k \\
        g_u
    \end{bmatrix} 
    \quad\text{and}\quad
    \matrixForm{B} = \begin{bmatrix}
        \matrixForm{B}_k \\
        \matrixForm{B}_u
    \end{bmatrix}.
\end{equation*}
We solve for the least squares solution
\begin{equation}\label{eq:exp_solve_for}
    \matrixForm{B}_k^* \cdot \matrixForm{B}_k \cdot \hat{c} = \matrixForm{B}_k^* \cdot g_k 
\end{equation}
and then
\begin{equation}\label{eq:exp_infer}
    g_u = \matrixForm{B}_u \cdot \hat{c}. 
\end{equation}
Here $\matrixForm{B}^*$ is the adjoint of $\matrixForm{B}$ whose action is to transform a sinogram $g$ to a vector whose components are the inner products of $g$ with the different basis functions $e^{ik\theta}U_n(s)W(s)$.

Due to the severe ill-conditioned nature of limited data tomography, the solution to \cref{eq:exp_solve_for} is not useful unless some kind of regularisation is used. We used Tikhonov regularisation, which means that enforcing the regularized range condition problem amounts to ridge regression with the basis functions from \eqref{eq:sino_synthesis}. As the orthonormal family of polynomials, we tried both the Legendre and Chebyshev polynomials and found that Chebyshev polynomials yielded results that were significantly more robust.

Notice that just as for the back projection operator we do not need to store the explicit matrix $\matrixForm{B}$, we only need to consider its action on coefficient vectors $c$ and the action of its adjoint on sinograms $g$. We do need the matrix $\matrixForm{B}_k^* \cdot \matrixForm{B}_k$, the size of which depends on our choice of $N$. We found that $N = 50$ yielded an adequate approximation which results in a manageable $650 \times 650$ matrix. Furthermore, with $N=50$ the operators $\matrixForm{B}$ and $\matrixForm{B}^*$ are of a computational complexity that is lower than the back projection operator resulting in a reconstruction procedure with a computational burden of the same magnitude as a normal \ac{FBP} reconstruction.

The computation of $\matrixForm{B}_k^* \cdot \matrixForm{B}_k$ is a relatively computationally expensive operation, however, this matrix only needs to be computed once during training and can then be stored for use during inference. Hence, computation is cheap compared to the computational load of training the network.

\begin{figure}
\begin{center}
  \begin{tikzpicture}[thick,scale=0.75, every node/.style={scale=0.75}]
    \input{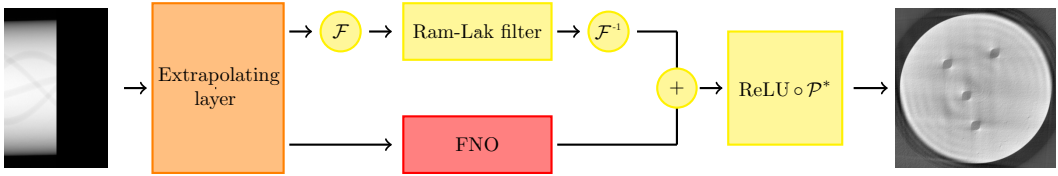}
  \end{tikzpicture}
\end{center}    
\caption{Architecture of the \ac{FNOBP} network}
\label{fnoarchitecture}
\end{figure}

\subsection{Network Architecture}
With no prior on the data distribution, enforcing range conditions is in some sense the best we can do. However, since the data comes from real-world objects, there is additional information hidden in the data, that we can exploit. The idea is to first make an initial extrapolation that preserves range conditions and then compensate for the deficiencies of the extrapolation with a neural network that learns from data. A straightforward generalization of the \ac{FBP} reconstruction operator is to replace convolution in data space with the reconstruction kernel $k$ with a convolution network, which can be adapted by training against a large training set. 

Since filters used in \ac{FBP} are usually easier to represent in the Fourier domain it seems plausible to use a network that operates by filtering in the frequency domain. An example of such an architecture that has been successful at other tasks is the \ac{FNO} \cite{fno_article}. The \ac{FNO} architecture was introduced in the context of operator learning. More specifically, the idea was to learn a neural network that outputs the solution of a partial differential equation from its initial/boundary conditions in a mesh-independent way (can be trained on one mesh and evaluated on another). In addition, the \ac{FNO} architecture is well-adapted for learning Fourier integral operators. Both the inverse of the ray transform and the \ac{FBP} are examples of Fourier integral operators \cite{krishnan2015microlocal}.

The key component of the \ac{FNO} is spectral convolution, which is essentially a usual convolutional layer but with every convolution performed in the frequency domain. Given an input with $C_{\text{in}}$ channels, regarded as a function $f \colon [0,1] \to \real^{C_{\text{in}}} \in L_1([0,1],\mathbb{R}^{C_{\text{in}}})$ the spectral convolutional layer yields an output $h \colon [0,1] \to \real^{C_{\text{out}}} \in L_1([0,1], \mathbb{R}^{C_{\text{out}}})$ of arbitrarily many channels $C_{\text{out}}$. For each channel $i$ the output is given by
\[
    \hat{h}(\omega)_i = \sum_{j = 1}^{C_{\text{in}}} 
\hat{k}_{ij}(\omega)\hat{f}(\omega)_j  
\]
which is equivalent to convolving the input $f$ with the filters $k_{ij} = \Fourier^{-1}(\hat{k}_{ij})$. Since the function domain is compact we consider discrete frequencies $\omega$.

The \ac{FNO} in our \ac{FNOBP} network treats every column of a sinogram, i.e. every set of measurements collected from the same rotation angle, as an individual channel. In this context, a sinogram sampled from $L$ directions is regarded as discrete samples from a function $\real \to \real^L$ that maps the spatial coordinate $u$ to the corresponding set of samples $g(:, u)$:
\[
u \mapsto \begin{bmatrix}
    g(\beta_1, u) \\
    \vdots  \\
    g(\beta_L, u)
\end{bmatrix}.
\]

In summary, the model can be written as
\[ g \mapsto \ReLU\Bigl( 
      \Radon^*\bigl(k_{\text{RL}} \ast g + \FNO(g)\bigr) 
    \Bigr) 
\]
where $k_{\text{RL}}$ is the commonly used Ram-Lak filter, see \cref{fnoarchitecture} for an illustration of the corresponding network architecture.
The \ac{FNO} layer itself consists of a set of spectral convolution layers with \ac{GELU} activation functions in between as well as linear skip connection layers, see \cite{fno_article}. After exploring a few configurations of hyperparameters we settled on using an \ac{FNO} with three hidden layers each consisting of $C = 60$ channels, each kernel $\hat{k}_{ij}$ was parameterized by $280$ values corresponding to the first $280$ Fourier coefficients of the filter $k_{ij}$.

The FNO exists in 1D and 2D variants, and although our data is inherently 2D (sinograms are functions on lines and angles), our FNO-BP employs the 1D variant for two reasons. First off, lower dimensional Fourier transforms are more efficient and secondly, the FNO performs translation invariant transforms, which means it would incorrectly treat all regions of the sinograms equally. This is not the case, since limited angle tomography results in missing or low quality data in concentrated regions that need special attention.

\subsection{Implementation}
The model was implemented in Python using PyTorch and ODL \cite{Odlgroup_2018}, and the code can be found on GitHub \cite{GitHubRepo}. The models were trained for 30 epochs on a synthesized dataset of 2500 phantoms, using the Adam optimizer with a learning rate of $3\cdot 10^{-5}$ and otherwise default PyTorch parameters. The synthesized data consisted of circular phantoms with different types of simple holes in them. Note that the phantom centers varied and the holes were of differing size and geometry, split between rectangular and oval. This is similar to the test data from HTC-2022, although the test data contains some phantoms with considerably more complex interiors.

\section{Results} \label{results}
The \ac{FNOBP} model was trained against synthetic data \cite{GitHubRepo} constructed to emulate the teaching data available from the HTC-2022 challenge \cite{HTC2022_data}. The \ac{FNOBP} model was then evaluated on the HTC-2022 testing data \cite{HTC2022_data}, together with a standard \ac{FBP} and an \ac{FBP} utilizing the sinogram extrapolation described in \cref{sec:sino_exp}. The resulting reconstructions had all negative pixels set to zero and were thresholded using Otsu's method \cite{Otsu_method}. Scoring, as displayed in \cref{result_table}, was done using Matthew's correlation coefficient, as done during the HTC-2022 \cite{HTC2022_rules} competition.

In \cref{result_table}, the scores for \ac{FNOBP} model can be seen in comparison to a regular \ac{FBP}, with and without sinogram extrapolation, as well as the winning submission to the HTC-2022 \cite{TUD_HHUD_HTC2022}. The methods were tested for different amounts of available angles (from $90^\circ$ to $30^\circ$), and the scores were taken as the mean over the three test images that were given for each difficulty level. It should be noted for both tables, that already a $90^\circ$ angular span (half of the required $180^\circ$) is quite limited. 

In \cref{fig:result_table_images:seg}, we see the resulting reconstructions for different data in the testing dataset. The images are ordered according to the method of reconstruction and available angular span. From these images, we can see that both methods utilizing the range conditions through sinogram extrapolation (\ac{FNOBP} and \acs{FBP} + range cond.) maintain the general shape of the disc. \ac{FNOBP} was trained exclusively on circular phantoms, explaining this fact, but \acs{FBP} + range cond. was not, meaning some part of the shape-preserving ability is possibly due to the sinogram extrapolation. Note that the centers of the phantoms differ between images and that this is accurately captured in the reconstructions. Both \ac{FBP} methods smear the reconstructions, however the smearing is contained within the disc when using range conditions. While the \ac{FBP} methods seem to exaggerate the holes for lower angle spans, \ac{FNOBP} instead seems to make them smaller. 

\begin{table}[tbh!]
\begin{center}
\begin{tabular}{l | c c c c c c c}
 Methods & $90^\circ$ & $80^\circ$ & $70^\circ$ & $60^\circ$ & $50^\circ$ & $40^\circ$ & $30^\circ$ \\
\hline
 \acs{FBP} & 0.654 & 0.685 & 0.634 & 0.614 & 0.520 & 0.394 & 0.284 \\
%\hline
 \acs{FBP} + range cond. & 0.851 & 0.797 & 0.689 & 0.667 & 0.612 & 0.492 & 0.404 \\
%\hline
 \acs{FNOBP} & 0.913 & 0.919 & 0.828 & 0.832 & 0.832 & 0.715 & 0.630 \\
%\hline
 Top HTC 2022 & 0.986 & 0.990 & 0.977 & 0.972 & 0.975 & 0.938 & 0.803 
\end{tabular}
\end{center}
\caption{Performance scores for reconstruction methods based on the metric used in the HTC~2022 competition. Scores are in $[-1,1]$ with 1 indicating a perfect reconstruction, 0 representing a reconstruction no better than random noise, and $-1$ a reconstruction that disagrees with the ground truth at every pixel.
Top HTC~2022 refers to the method that won the HTC~2022 competition.}
\label{result_table}
\end{table}

From \cref{result_table}, it seems that \ac{FNOBP} would achieve better reconstructions than either of the \ac{FBP} methods, but worse than the winning submission to HTC-2022 \cite{TUD_HHUD_HTC2022}. \Ac{FNOBP}'s capability compared to the \ac{FBP} methods is confirmed in \cref{fig:result_table_images,fig:result_table_images:seg}, where for the chosen ground truths \ac{FNOBP} achieves similarly accurate reconstructions to the winning submission from HTC-2022 down to an angle span of $50^\circ$, while both standard \ac{FBP} and sinogram extrapolation followed by \ac{FBP} begin to fall behind at angles spans of $80^\circ$ or $90^\circ$. For the smallest angles, all three of the above mentioned methods seem to perform poorly compared to the winning HTC-2022 submission.

The observations above are reinforced by \cref{fig:result_table_images,fig:result_table_images:seg}, where we study the different methods' ability to reconstruct a specific object, but for different angular spans. Here we also forego any thresholding. From these images, we confirm that even without thresholding the methods utilising range conditions largely maintain the shape of the object after reconstruction, in particular, \ac{FNOBP} maintains the shape very closely without much smearing outside of the object radius even for low angular spans. The discernibility of the holes is similar for all methods, although \ac{FBP} with range conditions seems to introduce more artefacts, in the shape of rings around the object's centre.

\section{Conclusion}\label{conclusion}
\Ac{FNOBP} is a novel neural network architecture for tomographic reconstruction that is designed to have a computational complexity during inference that is of the same magnitude as an \ac{FBP} reconstruction. 
Hence, \ac{FNOBP} is computationally much more efficient than model-based methods (iterative and variational) and contemporary domain informed learned methods such as LPD \cite{Adler2018} that require multiple forward- and back projections. 

\Ac{FNOBP} was tested on limited angle tomography data in the HTC-2022 challenge \cite{HTC2022_data}. While \ac{FNOBP} is outperformed by the winning submission in HTC-2022 \cite{TUD_HHUD_HTC2022}, the winning submission makes use of a much larger network to learn a direct mapping between the set of sinograms and reconstructions. We believe that \ac{FNOBP} has the potential for better generalisation properties and requires less data and computational power to train.
Furthermore, although we used \ac{FNOBP} for limited angle reconstruction, this reconstruction network is also applicable to other forms of limited data tomography, like truncated or sparse view data.

Directions for further explorations could include using \ac{FNO} in 3D settings or exploiting its discretization invariance to learn reconstruction operators that are transferable between different environments. Furthermore, it would be interesting to explore further how neural networks for sinogram inpainting could be combined with analytical extrapolation procedures that make use of range conditions.

\printbibliography[title={References},heading=bibintoc]

\clearpage

\begin{table}[tbh!]
\centering
\begin{tabular}{ccccc}
Ground truth & & \acs{FNO} & \acs{FBP} + range cond. & \acs{FBP}  
\\
\imagetop{\includegraphics[width=\Rwidth\linewidth]{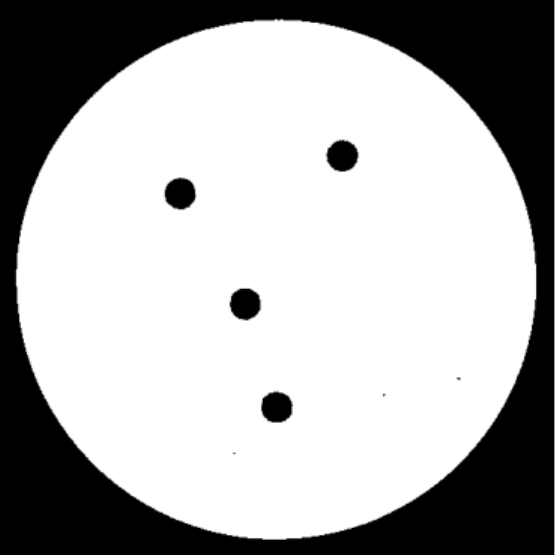}} &
\imagetop{$90^\circ$} & 
\imagetop{\includegraphics[width=\Rwidth\linewidth]{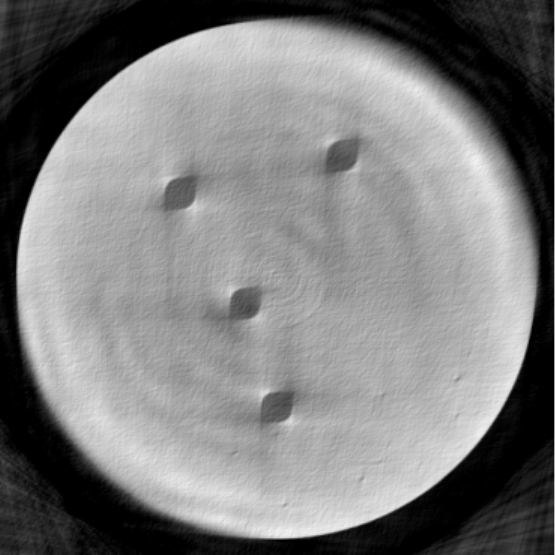}} &
\imagetop{\includegraphics[width=\Rwidth\linewidth]{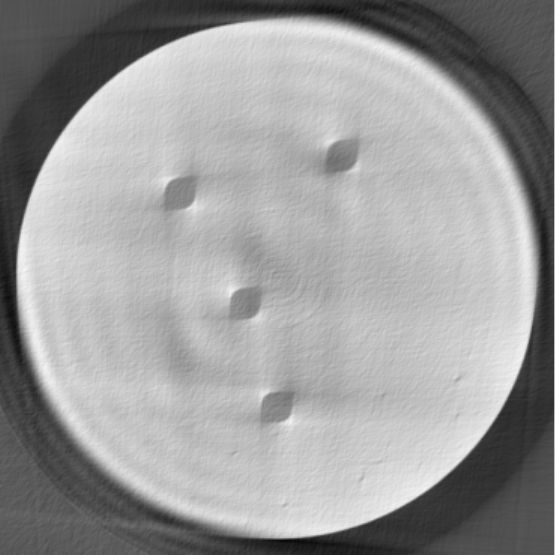}} &
\imagetop{\includegraphics[width=\Rwidth\linewidth]{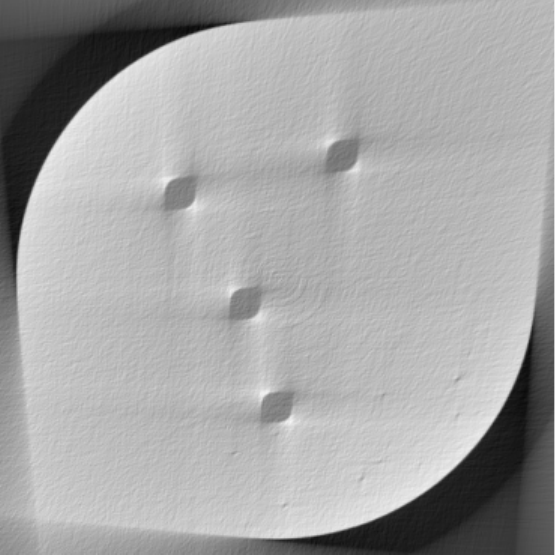}}
\\
&\imagetop{$80^\circ$} &
\imagetop{\includegraphics[width=\Rwidth\linewidth]{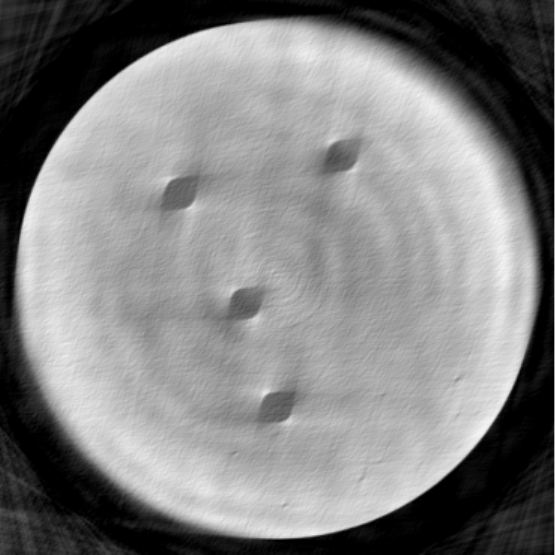}} &
\imagetop{\includegraphics[width=\Rwidth\linewidth]{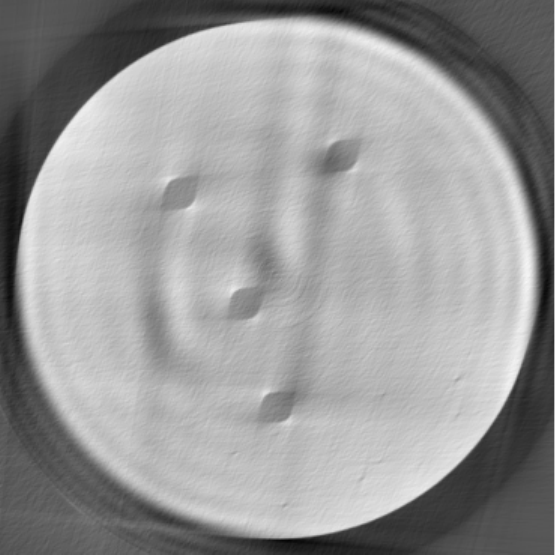}} &
\imagetop{\includegraphics[width=\Rwidth\linewidth]{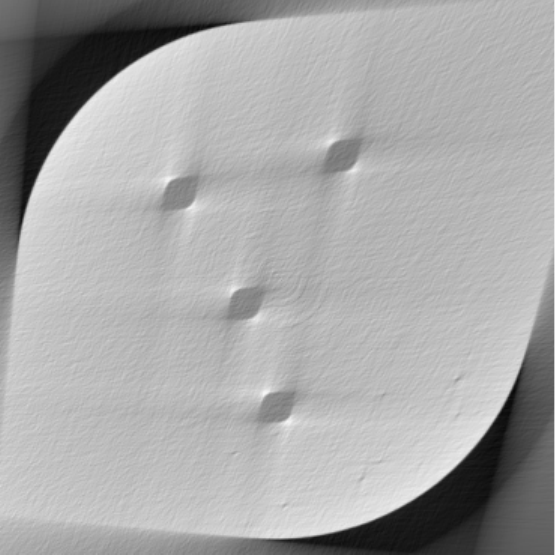}} 
\\
&\imagetop{$70^\circ$} &
\imagetop{\includegraphics[width=\Rwidth\linewidth]{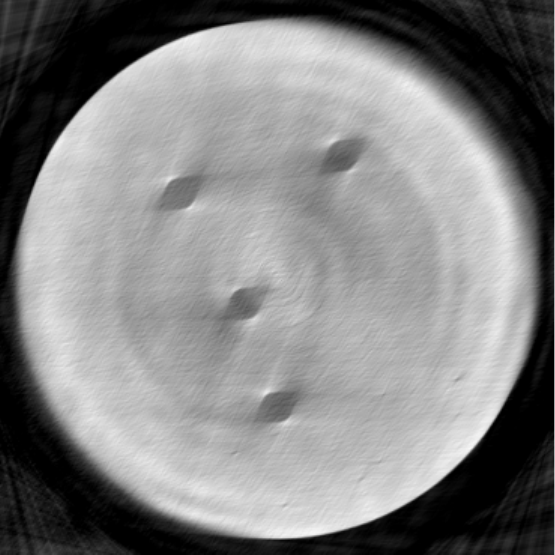}} &
\imagetop{\includegraphics[width=\Rwidth\linewidth]{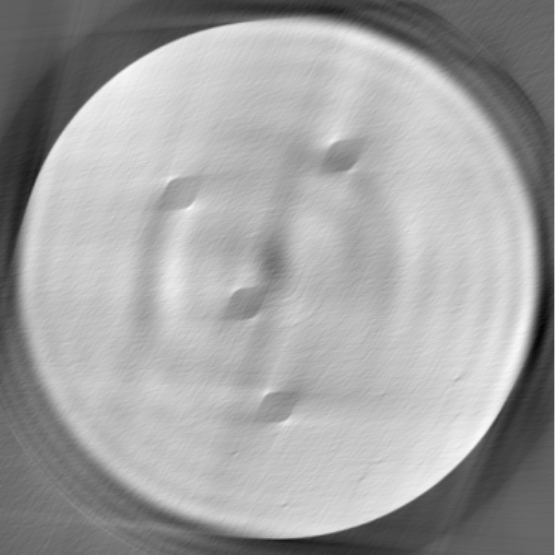}} &
\imagetop{\includegraphics[width=\Rwidth\linewidth]{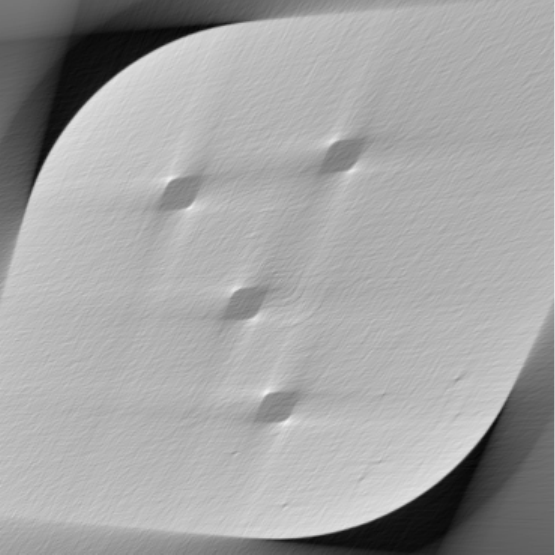}}
\\
&\imagetop{$60^\circ$} &
\imagetop{\includegraphics[width=\Rwidth\linewidth]{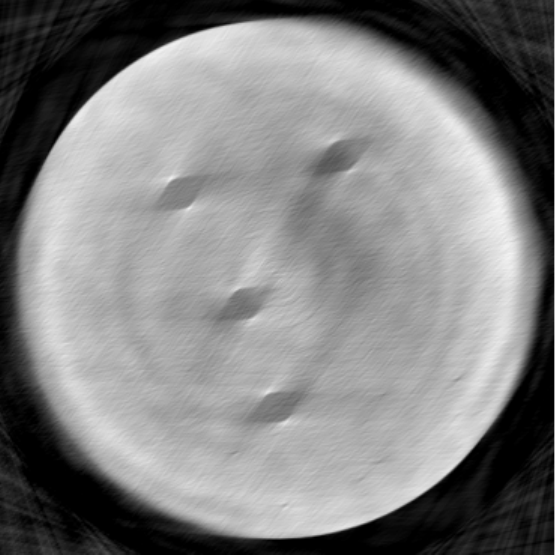}} &
\imagetop{\includegraphics[width=\Rwidth\linewidth]{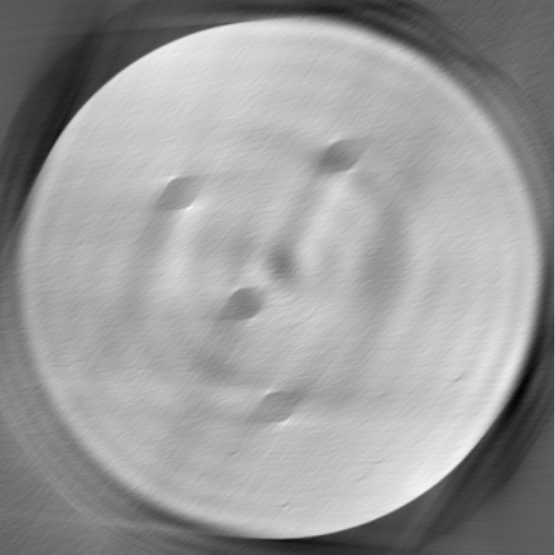}} &
\imagetop{\includegraphics[width=\Rwidth\linewidth]{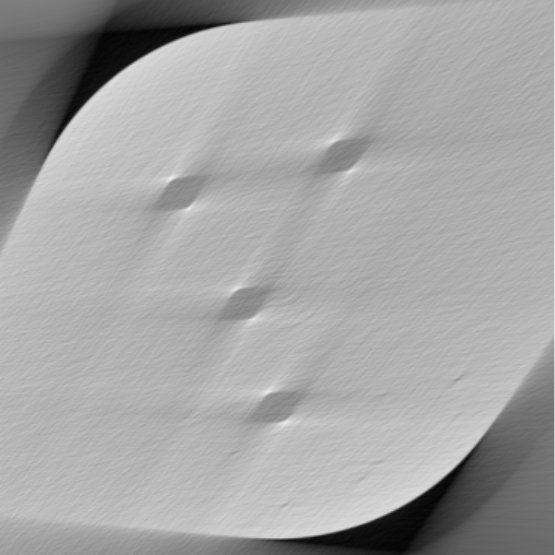}}
\\
&\imagetop{$50^\circ$} &
\imagetop{\includegraphics[width=\Rwidth\linewidth]{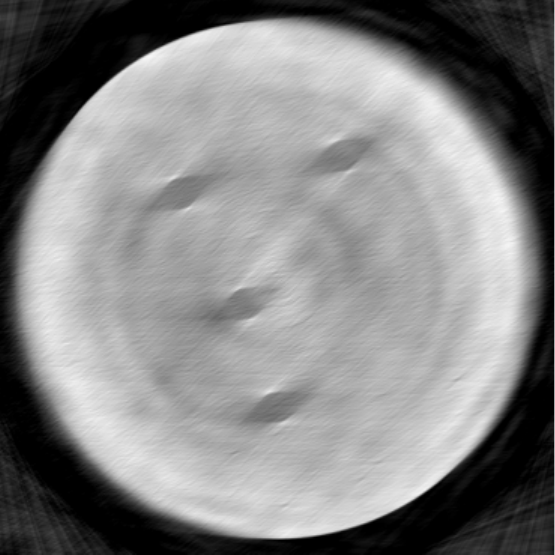}} &
\imagetop{\includegraphics[width=\Rwidth\linewidth]{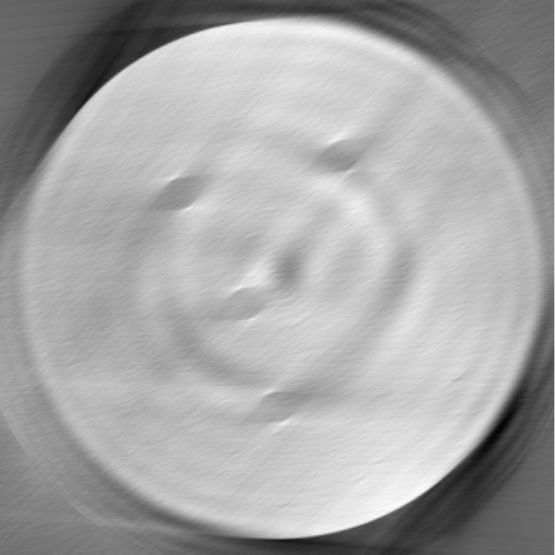}} &
\imagetop{\includegraphics[width=\Rwidth\linewidth]{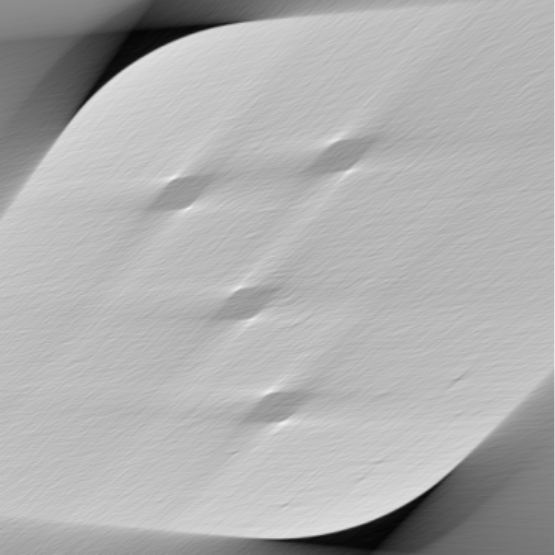}}
\\
&\imagetop{$40^\circ$} &
\imagetop{\includegraphics[width=\Rwidth\linewidth]{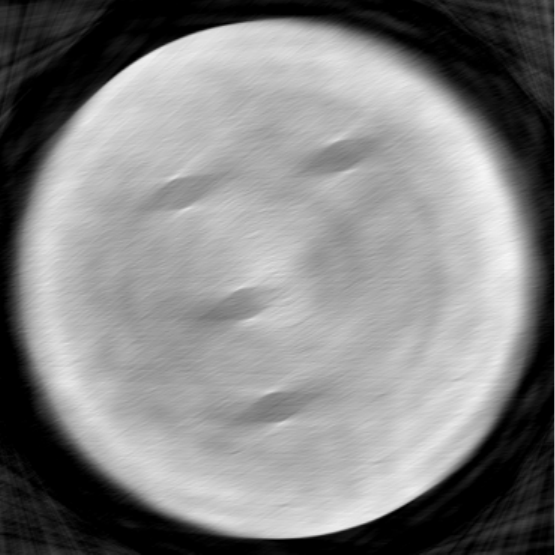}} &
\imagetop{\includegraphics[width=\Rwidth\linewidth]{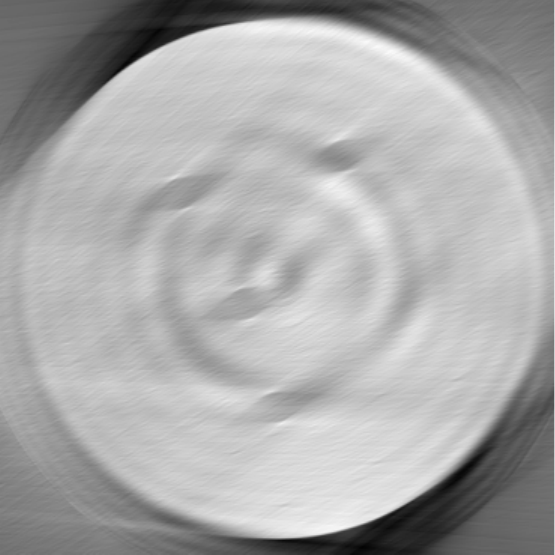}} &
\imagetop{\includegraphics[width=\Rwidth\linewidth]{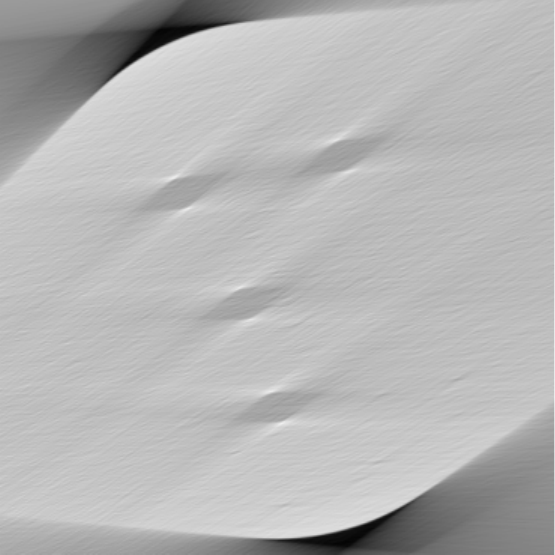}}  
\\
&\imagetop{$30^\circ$} &
\imagetop{\includegraphics[width=\Rwidth\linewidth]{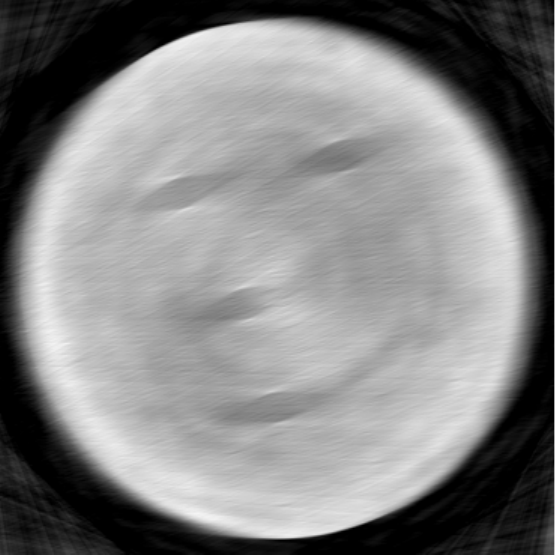}} &
\imagetop{\includegraphics[width=\Rwidth\linewidth]{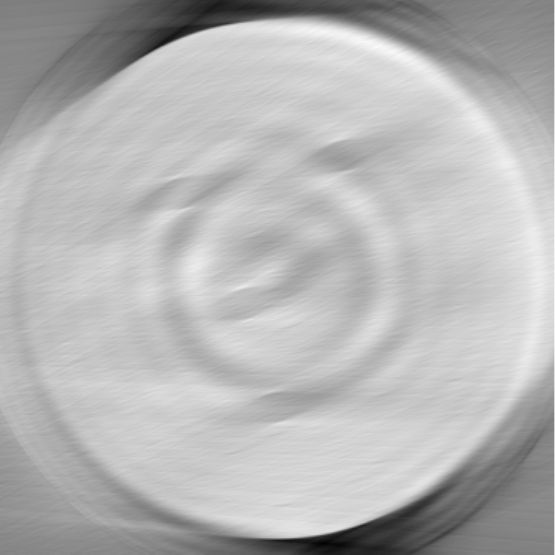}} &
\imagetop{\includegraphics[width=\Rwidth\linewidth]{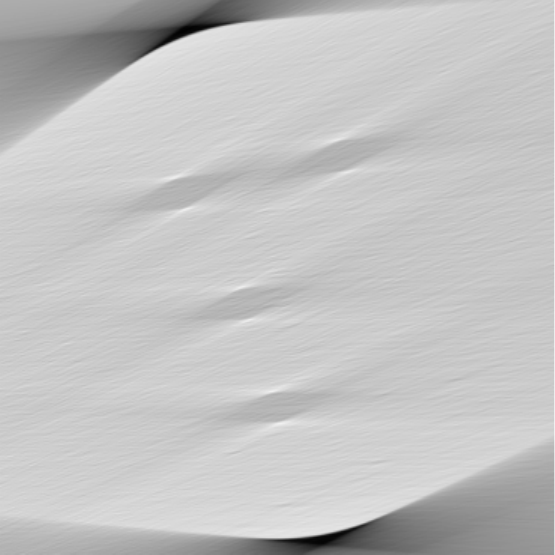}}
\end{tabular}
\caption{Reconstructions of an example phantom. The available angle span is displayed to the left ranging from $90^\circ$ (half of full data) to $30^\circ$. Ground truth is displayed at the top left.}
\label{fig:result_table_images}
\end{table}

\clearpage
\begin{table}[tbh!]
\centering
\begin{tabular}{ccccc}
Ground truth && \acs{FNO} & \acs{FBP} + range cond. & \acs{FBP}  
\\
\imagetop{\includegraphics[width=\Rwidth\linewidth]{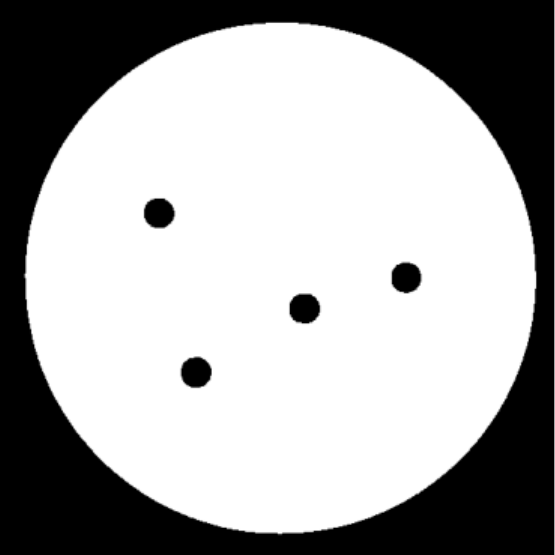}} &
\imagetop{$90^\circ$} & 
\imagetop{\includegraphics[width=\Rwidth\linewidth]{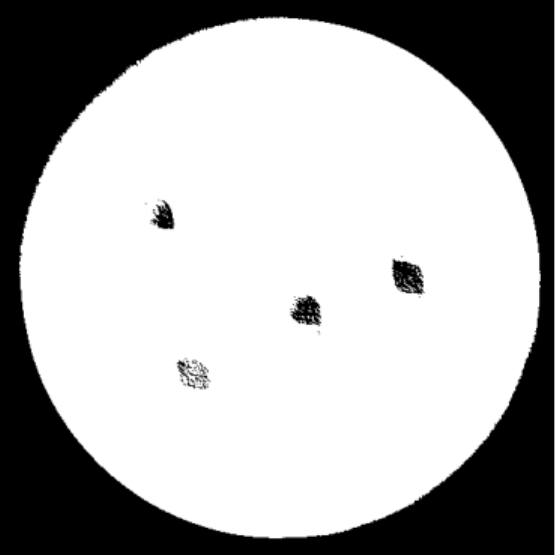}} &
\imagetop{\includegraphics[width=\Rwidth\linewidth]{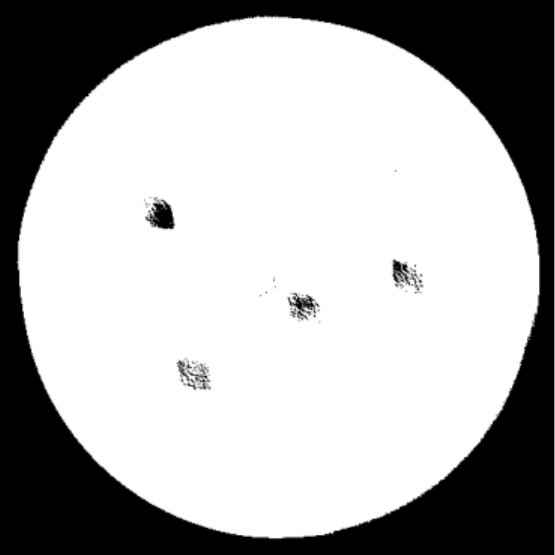}} &
\imagetop{\includegraphics[width=\Rwidth\linewidth]{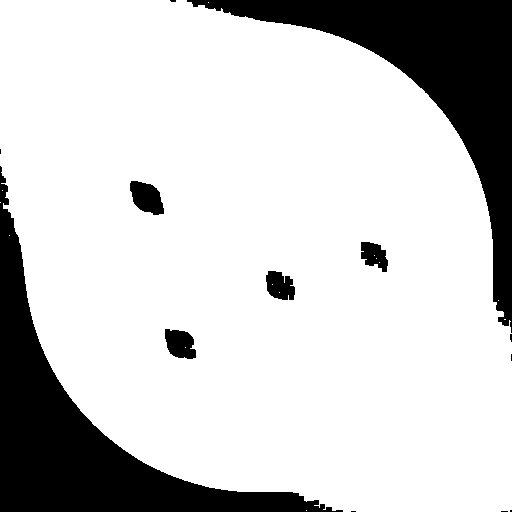}} 
\\
\imagetop{\includegraphics[width=\Rwidth\linewidth]{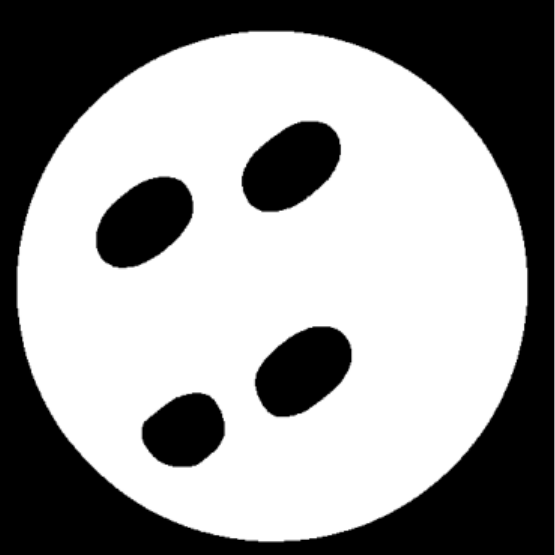}} &
\imagetop{$80^\circ$} &
\imagetop{\includegraphics[width=\Rwidth\linewidth]{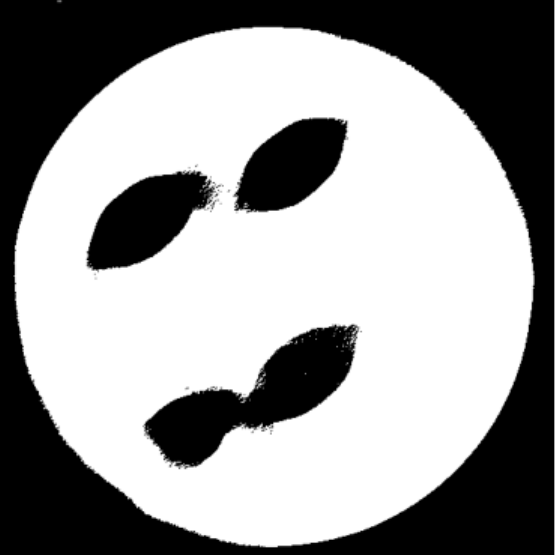}} &
\imagetop{\includegraphics[width=\Rwidth\linewidth]{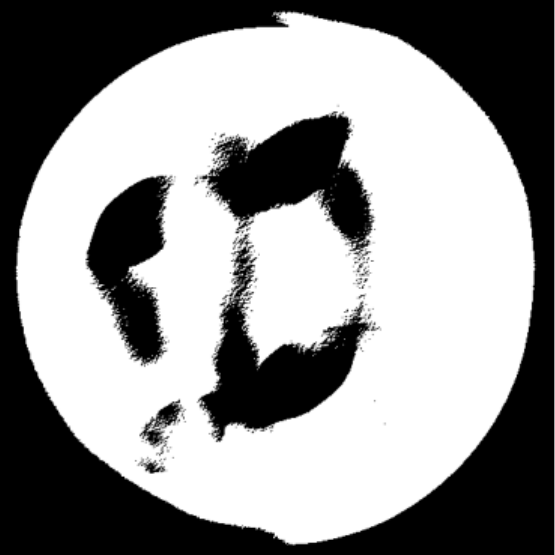}} &
\imagetop{\includegraphics[width=\Rwidth\linewidth]{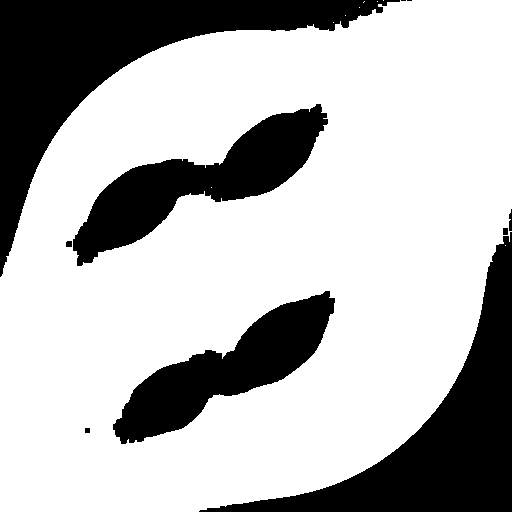}}  
\\
\imagetop{\includegraphics[width=\Rwidth\linewidth]{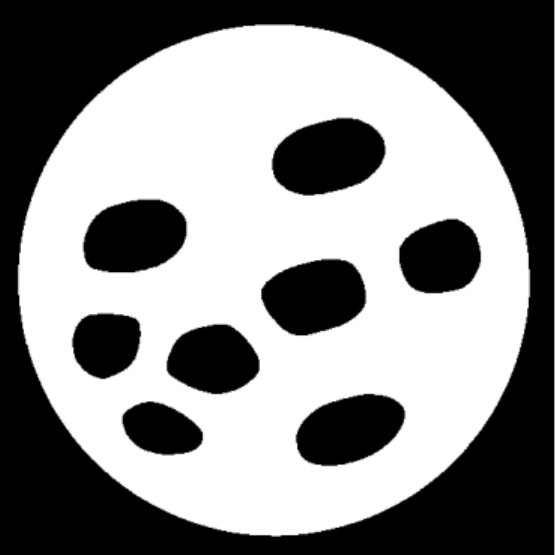}} &
\imagetop{$70^\circ$} &
\imagetop{\includegraphics[width=\Rwidth\linewidth]{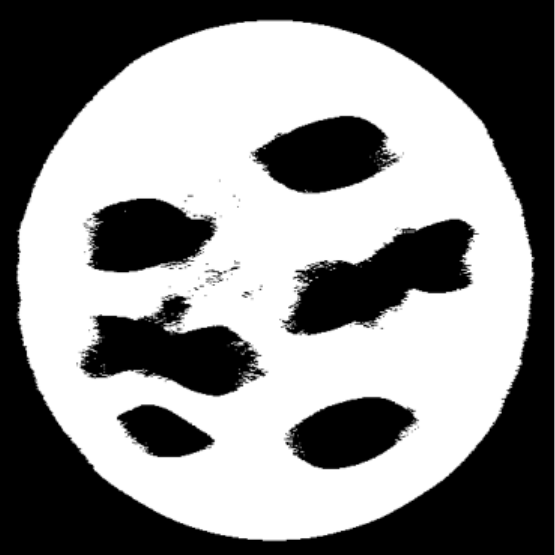}} &
\imagetop{\includegraphics[width=\Rwidth\linewidth]{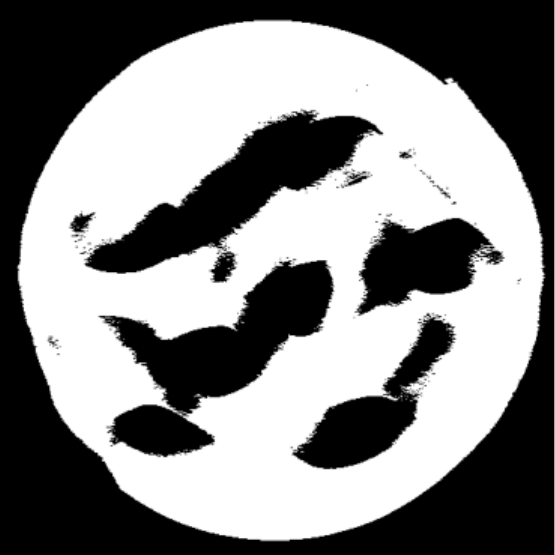}} &
\imagetop{\includegraphics[width=\Rwidth\linewidth]{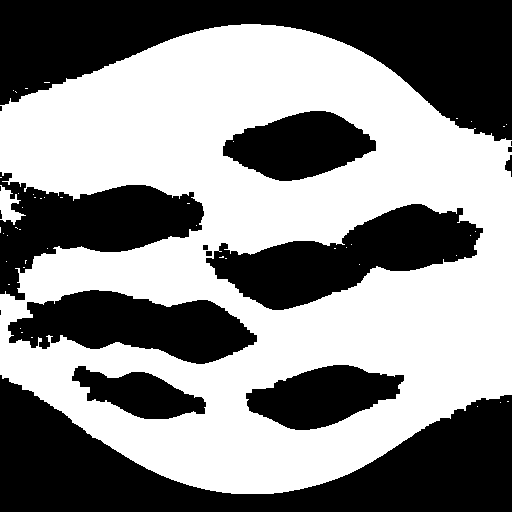}}
\\
\imagetop{\includegraphics[width=\Rwidth\linewidth]{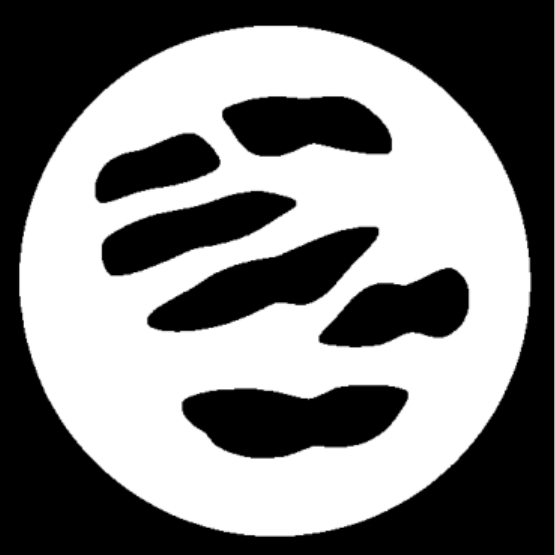}} &
\imagetop{$60^\circ$} &
\imagetop{\includegraphics[width=\Rwidth\linewidth]{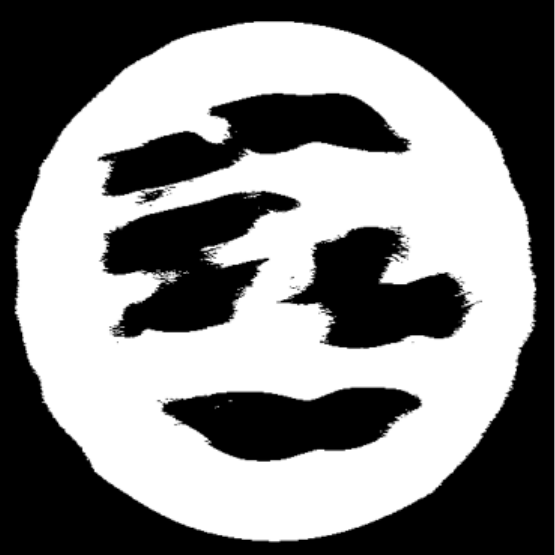}} &
\imagetop{\includegraphics[width=\Rwidth\linewidth]{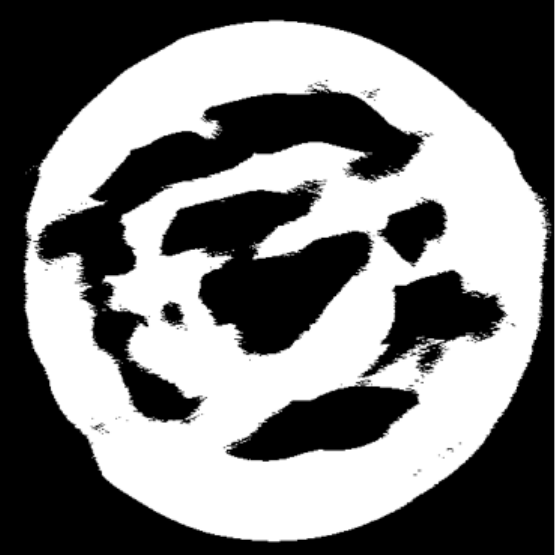}} &
\imagetop{\includegraphics[width=\Rwidth\linewidth]{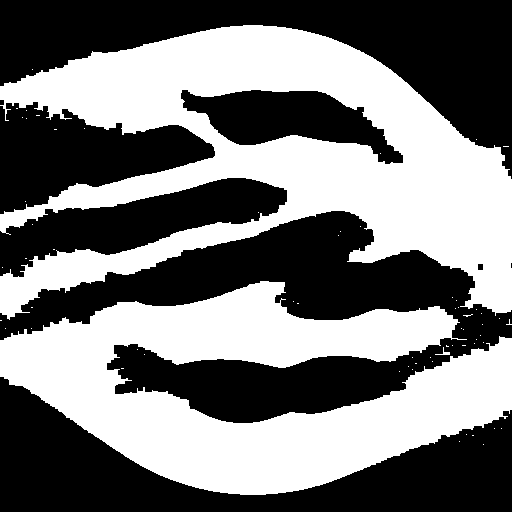}}
\\
\imagetop{\includegraphics[width=\Rwidth\linewidth]{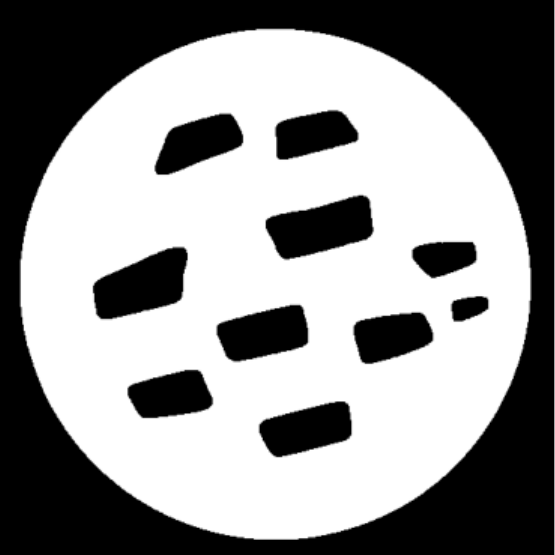}} &
\imagetop{$50^\circ$} &
\imagetop{\includegraphics[width=\Rwidth\linewidth]{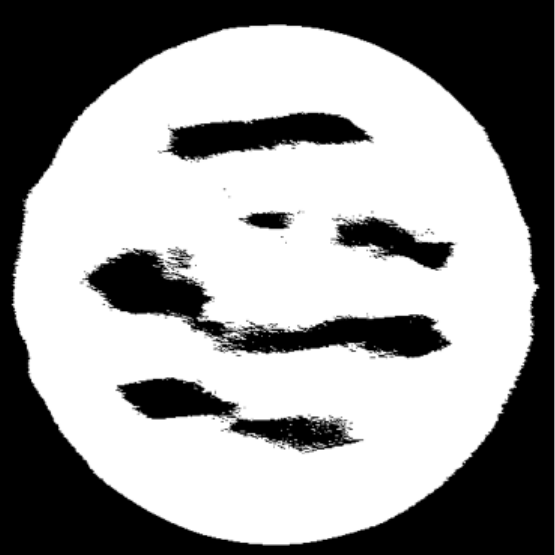}} &
\imagetop{\includegraphics[width=\Rwidth\linewidth]{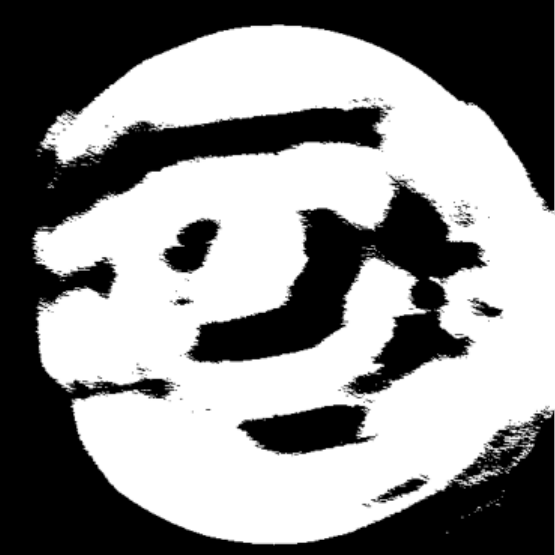}} &
\imagetop{\includegraphics[width=\Rwidth\linewidth]{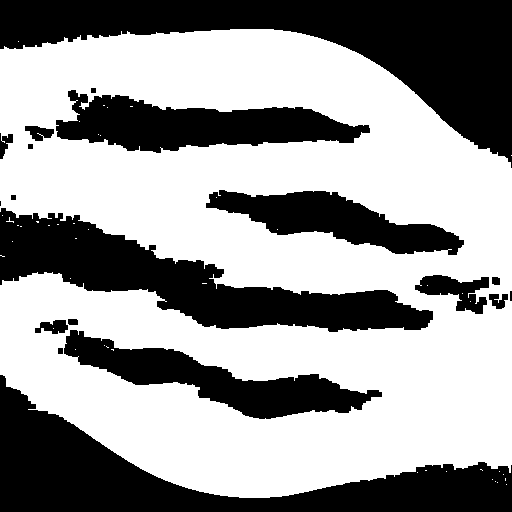}}
\\
\imagetop{\includegraphics[width=\Rwidth\linewidth]{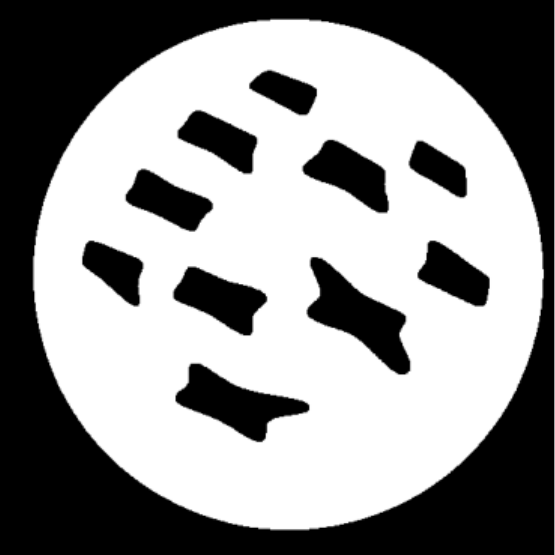}} &
\imagetop{$40^\circ$} &
\imagetop{\includegraphics[width=\Rwidth\linewidth]{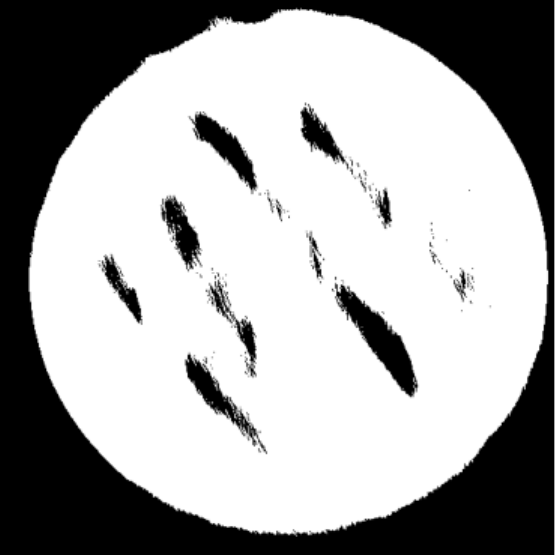}} &
\imagetop{\includegraphics[width=\Rwidth\linewidth]{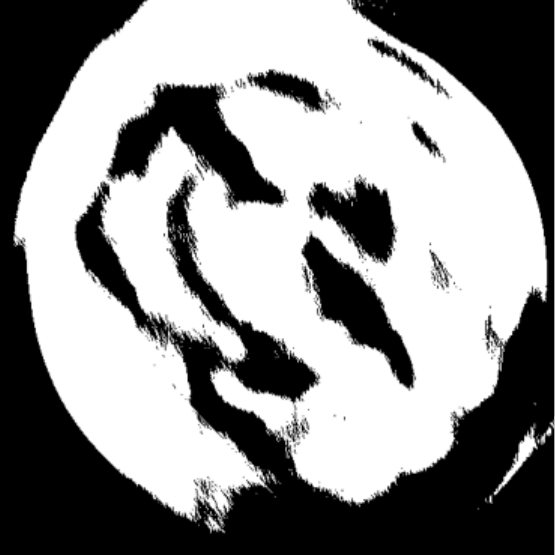}} &
\imagetop{\includegraphics[width=\Rwidth\linewidth]{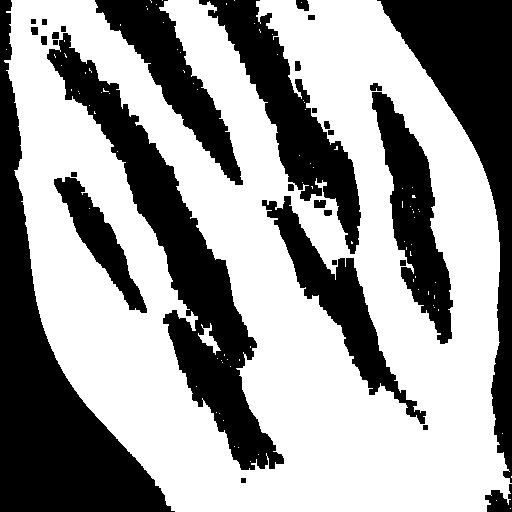}}
\\
\imagetop{\includegraphics[width=\Rwidth\linewidth]{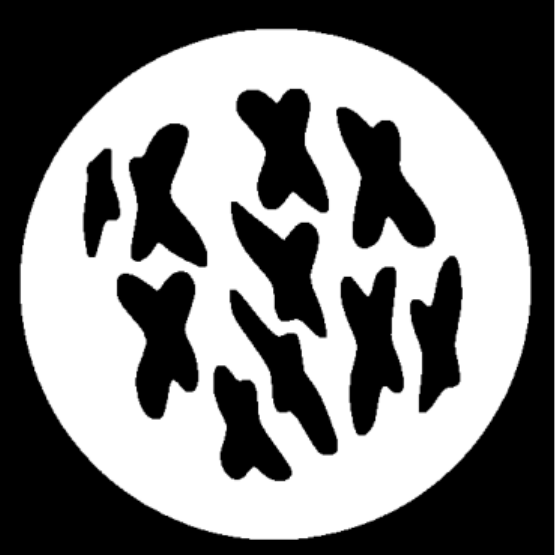}} &
\imagetop{$30^\circ$} &
\imagetop{\includegraphics[width=\Rwidth\linewidth]{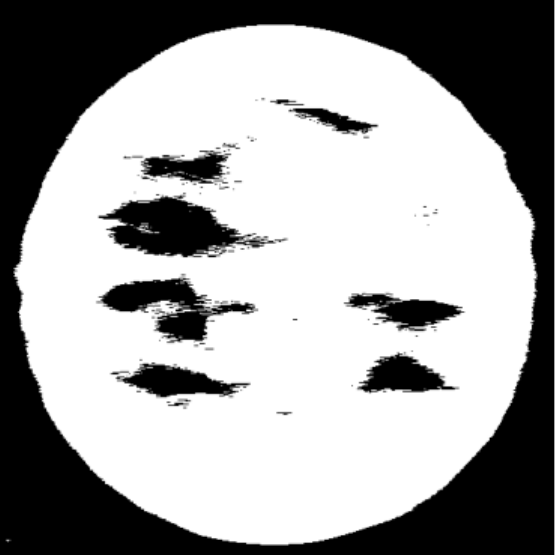}} &
\imagetop{\includegraphics[width=\Rwidth\linewidth]{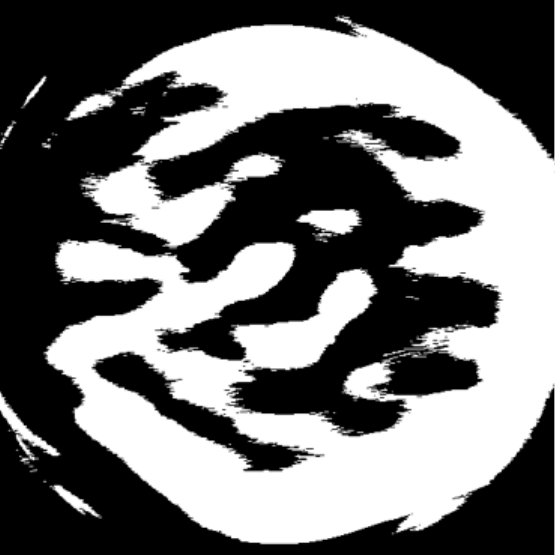}} &
\imagetop{\includegraphics[width=\Rwidth\linewidth]{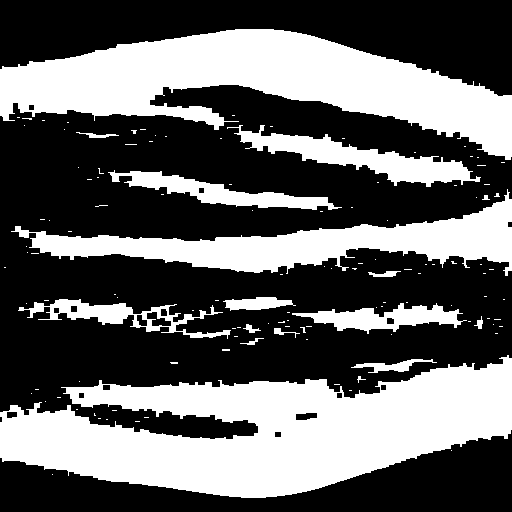}}
\end{tabular}
\caption{Reconstructions arranged according to increasing difficulty. Bottom row is most complicated object and largest missing data.}
\label{fig:result_table_images:seg}
\end{table}

\end{document}